\documentclass[11pt,A4]{amsart}

\usepackage[top=3cm,bottom=3cm,left=3cm,right=3cm]{geometry}
\usepackage[dvipsnames,svgnames,x11names]{xcolor}
\usepackage{amsmath,amssymb,mathtools}
\usepackage{graphicx}
\usepackage{float}
\usepackage[all]{xy}
\usepackage{verbatim}

\usepackage{color}

\newcommand{\hs}{\kern 0.8pt}
\newcommand{\hssh}{\kern 1.2pt}
\newcommand{\hshs}{\kern 1.6pt}
\newcommand{\hssss}{\kern 2.0pt}
\newcommand{\ha}{{\kern 1pt}}

\newcommand{\hm}{\kern -0.8pt}
\newcommand{\hmm}{\kern -1.2pt}

\newcommand{\upsig}{{\hs\sigma\!}}
\newcommand{\upg}{{\hs g\!}}

\begin{document}

\title[Non-abelian  group cohomology with coefficients in a crossed module]
{Non-abelian hypercohomology of a group\\ 
with coefficients in a crossed module,\\ 
and Galois cohomology}

\address{School of Mathematical Sciences, Tel Aviv University, 6997801 Tel Aviv, Israel}
\email{borovoi@tauex.tau.ac.il}.
\author{Mikhail Borovoi}

\thanks{Partially supported by NSF Grant DMS-8610730}

\begin{abstract}
We develop a hypercohomology theory of a group with coefficients in a crossed
module, and apply it to define abelianization maps for Galois cohomology of reductive
algebraic groups.

This is an arXiv version of my old preprint (Princeton 1991 -- Bonn 1992).
Since then, Sections 1--3 were covered in the paper of B. Noohi [Noo] written at my request,
and Section 4 was covered in my paper [Bo3].
\end{abstract}

\maketitle

\section*{Introduction}

Let $\Gamma$ be a profinite group and let
\begin{equation*}
1 \to \stackrel{\scriptstyle -1}{F} \stackrel{\scriptstyle \alpha}{\longrightarrow} \stackrel{\scriptstyle 0}{G} \to 1
\end{equation*}
be a short complex (a complex of length 2) of (in general non-abelian) groups, where the
numbers $-1$ and $0$ over the letters denote the degrees: $F$ is in degree $-1$ and $G$ is in
degree $0$. We assume that the group $\Gamma$ acts on $F$ and $G$, and that $\alpha$ is a homomorphism
of $\Gamma$-groups.

For applications to the Galois cohomology of connected algebraic groups, we would
like to be able to define the first hypercohomology set $\mathbf{H}^{1}(\Gamma, F \to G)$ in a functorial way.
In general this is not likely to be possible. Indeed, if we take $G$ to be $\{1\}$, then we must
have $\mathbf{H}^{1}(\Gamma, F \to 1) = H^{2}(\Gamma, F)$. However, as far as I know, there is no functorial definition
of second cohomology in the non-abelian case (the second nonabelian cohomology theories of Springer
[Sp] and Giraud [Gi] are evidently non-functorial).

Fortunately, it is possible to define the first hypercohomology in a functorial way when
$F \to G$ is \textit{a crossed module}. A crossed module is a group homomorphism $F \stackrel{\scriptstyle\alpha}{\longrightarrow} G$ with
an action of $G$ on $F$ satisfying certain natural conditions (see 2.1 for the precise definition
and [BHu] for a survey). The notion of a crossed module was introduced in 1946 by
J.H.C.Whitehead [W1], [W2], who was motivated by topological problems.

\indent To write down hypercohomology exact sequences, what we need is not only $\mathbf{H}^{1}$, but
also $\mathbf{H}^{-1}$ and $\mathbf{H}^{0}$. In Section 1 for any short complex of $\Gamma$-groups $F \to G$ (not necessarily
a crossed module), we define, in terms of cocycles, an abelian group $\mathbf{H}^{-1}(F \to G)$ and
pointed set $\mathbf{H}^{0}(F \to G)$, where we write $\mathbf{H}^{i}(F \to G)$ for $\mathbf{H}^{i}(\Gamma, F \to G)$. These definitions
were earlier given by Deligne [Del], 2.4.3, in terms of torsors. (Deligne writes $\mathbf{H}^{0}$ for our
$\mathbf{H}^{-1}$, and $\mathbf{H}^{1}$ for our $\mathbf{H}^{0}$.)

For a crossed module $F \to G$ with $\Gamma$-action, we define in Section 2 a group structure
on $\mathbf{H}^{0}(\Gamma, F \to G)$. Then we define, again in cocyclic form, the first hypercohomology
set $\mathbf{H}^{1}(\Gamma, F \to G)$. We follow Dedecker [Ded2],[Ded3], who defined $\mathbf{H}^{1}(\Gamma, F \to G)$ for a
crossed module $F \to G$ with trivial $\Gamma$-action; the generalization to the case of non-trivial
$\Gamma$-action is obvious. Note that Dedecker regards $\mathbf{H}^{1}(\Gamma, F \to G)$ not as hypercohomology
of a complex, but as a nice, functorial definition of the second cohomology $H^{2}(\Gamma, F)$,
so the group $G$ and its action on $F$ are for him just auxiliary structures necessary to
define $H^{2}(\Gamma, F)$. Dedecker writes $H^{2}(\Gamma, F \to G)$ for our $\mathbf{H}^{1}(\Gamma, F \to G)$. We regard
$\mathbf{H}^{1}(\Gamma, F \to G)$ as hypercohomology of a complex, and write down the hypercohomology
exact sequence assigned to  a short exact sequence of crossed modules.

In Section 3 we use hypercohomology exact sequences to prove

\medskip

\noindent\textbf{Theorem} (Theorem 3.3) \textit{Let $(F_{1} \to G_{1}) \to (F_{2} \to G_{2})$ be a quasi-isomorphism of crossed
modules with $\Gamma$-action. Then the induced maps $\mathbf{H}^{i}(F_{1} \to G_{1}) \to \mathbf{H}^{i}(F_{2} \to G_{2})$ $(i =$
$-1, 0, 1)$ are bijections.}

\medskip

In Section 4 we apply results of Sections 1--3 to the crossed module of algebraic groups
$G^{\mathrm{sc}} \overset{\rho}{\longrightarrow} G$, introduced by Deligne ([Del], 2.4.7). Here G is a connected reductive algebraic
group over a field $K$ of characteristic $0$, $G^{\mathrm{sc}}$ is the universal covering of the derived group
$G^{\mathrm{ss}}$ of $G$, the homomorphism $\rho$ is the composition $G^{\mathrm{sc}} \to G^{\mathrm{ss}} \to G$, and $G$ acts on $G^{\mathrm{sc}}$ in
the obvious way. Let $Z$ be the center of $G$ and $Z^{(\mathrm{sc})}$ the center of $G^{\mathrm{sc}}$. Let $H^{0}(K, G)$ and
$H^{1}(K, G)$ denote the $0$-dimensional and $1$-dimensional Galois cohomology of $G$. We define
the abelian Galois cohomology groups of $G$ by
\begin{equation*}
H^{i}_{\mathrm{ab}}(K, G) = \mathbf{H}^{i}(K, Z^{(\mathrm{sc})} \to Z) \qquad (i \geq -1).
\end{equation*}
Using the morphism $(1 \to G) \to (G^{\mathrm{sc}} \to G)$ and the quasi-isomorphism $(Z^{(\mathrm{sc})} \to Z) \to$
$(G^{\mathrm{sc}} \to G)$ of crossed modules of algebraic groups, we define for $i = 0, 1$ the abelianization maps
\begin{equation*}
\mathrm{ab}^{i}\colon H^{i}(K, G) \to \mathbf{H}^{i}(G^{\mathrm{sc}} \to G) \overset{\sim}{\longrightarrow} \mathbf{H}^{i}(Z^{(\mathrm{sc})} \to Z) = H^{i}_{\mathrm{ab}}(K, G).
\end{equation*}
The abelianization map $\mathrm{ab}^{0}$ was first defined by Deligne [Del]. The map $\mathrm{ab}^{1}$ generalizes
a map of Kottwitz ([Ko2], Thm. 1.2), which he defined and extensively used in the case
when $K$ is a local field. Kottwitz defined the abelianization map with the help of a
rather complicated method of $z$-extensions of reductive groups. The hypercohomology
with coefficients in a crossed module permits us to define the maps $\mathrm{ab}^{0}$ and $\mathrm{ab}^{1}$ explicitly,
in particular in terms of cocycles (Propositions 4.3.1 and 4.3.2).

Constructions of Section 4 are used in our  [Bo3] (cf. also [Bo1]),
where we describe ``explicitly'' the first Galois cohomology of a connected reductive group over a number field. Such constructions are useful in cohomological calculations related to Shimura varieties; cf. [Mi].

Note that it is also possible to define the abelianization map $\mathrm{ab}^{2}\colon H^{2}(K, G) \to H^{2}_{\mathrm{ab}}(K, G)$
(cf. [Bo2]), where $H^{2}(K, G)$ is the second non-abelian Galois cohomology set of Springer
[Sp] and Giraud [Gi]. If $K$ is a local field or a number field and $\eta \in H^{2}(K, G)$, then
$\mathrm{ab}^{2}(\eta) = 0$ if and only if $\eta$ is a neutral class, i.e., it corresponds to a split extension.

\medskip

\noindent\textit{Remarks} \quad (1) The cocyclic constructions of Sections 1--3 go through in the more general
case of hypercohomology of a simplicial set with coefficients in a family of crossed modules.

(2) In [Br1] Breen defines $\mathbf{H}^{-1}$, $\mathbf{H}^{0}$ and $\mathbf{H}^{1}$ in a uniform way for a sheaf of crossed
modules $F \to G$ on a site, and constructs the hypercohomology exact sequence (0.1.1).
Breen uses the machinery of homotopical algebra. In the particular case of the site of $\Gamma$-sets
his definitions appear to be equivalent to ours. Our results were obtained independently
(before the paper [Br1] appeared).

(3) Breen ([Br1], 6.2) proves that $\mathbf{H}^{1}(F \to G)$ can be identified with the set of
equivalence classes of torsors under the Picard category assigned to the crossed module
$F \to G$. Deligne noticed (private communication) that our Theorem 3.3 follows from this
description of $\mathbf{H}^{1}(F \to G)$, because quasi-isomorphic crossed modules define equivalent
Picard categories.

(4) We claim no originality. Most of the results of Sections 1--3 are known (except the
construction of the connecting map $H^{1} \to H^{2}$ in 2.17 -- 2.22), see Remarks (2--3) above. We
need however this cocyclic exposition for Section 4, where we write down explicit cocyclic
formulas for $\mathrm{ab}^{0}$ and $\mathrm{ab}^{1}$.

\section*{1 \quad Hypercohomology in degrees $-1$ and $0$}

\noindent\textbf{1.1} \textit{Short complexes of groups.} \quad Let $\Gamma$ be a profinite group. A discrete $\Gamma$-group is a group
$G$ endowed with a left action of $\Gamma$ which is continuous with respect to the discrete topology
on $G$. Here ``continuous'' means that the stabilizer of any element $g \in G$ is open in $\Gamma$.
From now on, by a $\Gamma$-group we mean a discrete $\Gamma$-group.

Let $\alpha\colon F \to G$ be a morphism of $\Gamma$-groups, i.e., a group homomorphism respecting the
action of $\Gamma$. We consider $F \to G$ as a \textit{short complex}
\begin{equation*}
1 \to \stackrel{\scriptstyle -1}{F} \to \stackrel{\scriptstyle 0}{G} \to 1
\end{equation*}
where $F$ is in degree $-1$ and $G$ is in degree $0$.

\medskip

\noindent\textbf{1.2} \textit{Hypercohomology.} \quad We define hypercohomology in degree $-1$. We set
\begin{equation*}
\mathbf{H}^{-1}(F \stackrel{\alpha}{\longrightarrow} G) = (\ker\ \alpha)^{\Gamma}
\end{equation*}
where $( \ )^{\Gamma}$ means (the group of) invariants.

We define 0-hypercohomology. We write $\operatorname{Maps}(\Gamma, F)$ for the set of continuous maps $\varphi\colon \Gamma \to F$ and set
\begin{align*}
C^0 &= \mathrm{Maps}(\Gamma, F) \times G \text{ (we regard } C^0 \text{ as a set)}\\
Z^0 &= \{(\varphi, g) \in C^0 \mid \varphi(\sigma\tau) = \varphi(\sigma) \cdot {}^{\upsig}\varphi(\tau),\ {}^{\upsig}g = \alpha(\varphi(\sigma)^{-1}) \cdot g,\ \sigma, \tau \in \Gamma\}
\end{align*}
The group $F$ acts on the set of 0-cocycles $Z^0$ on the right by
\begin{equation*}
(\varphi, g) * f = (\varphi', g'), \quad \varphi'(\sigma) = f^{-1} \cdot \varphi(\sigma) \cdot {}^{\upsig}f, \quad g' = \alpha(f)^{-1} \cdot g,
\end{equation*}
(where $f \in F$), and we set
\begin{equation*}
\mathbf{H}^0(F \to G) = Z^0/F
\end{equation*}
The set $\mathbf{H}^0(F \to G)$ has a neutral element, namely the class of $(1,1)$. We write $\operatorname{Cl}(\varphi, g)$ for the hypercohomology class of a cocycle $(\varphi, g)$.

\textbf{1.3} \textit{Morphisms of complexes.} \quad A morphism of (short) complexes $(F_1 \to G_1) \to (F_2 \to G_2)$ is a commutative diagram
\begin{equation*}
\xymatrix{
F_1 \ar[r] \ar[d] & F_2 \ar[d] \\
G_1 \ar[r] & G_2
}
\end{equation*}
of $\Gamma$-groups. Such a morphism induces a canonical homomorphism
\begin{equation*}
\mathbf{H}^{-1}(F_1 \to G_1) \to \mathbf{H}^{-1}(F_2 \to G_2)
\end{equation*}
and a canonical map
\begin{equation*}
\mathbf{H}^0(F_1 \to G_1) \to \mathbf{H}^0(F_2 \to G_2).
\end{equation*}

\subsection*{1.4 Examples}

\begin{enumerate}
\renewcommand{\labelenumi}{(\arabic{enumi})}
\item $\mathbf{H}^0(1 \to G) = H^0(G) = G^\Gamma$.
\item $\mathbf{H}^0(F \to 1) = H^1(F)$. To $\operatorname{Cl}(\varphi, 1) \in \mathbf{H}^0(F \to 1)$ we assign $\operatorname{Cl}(\varphi) \in H^1(F)$.
\item If $\alpha\colon F \to G$ is injective, then the morphism of complexes $(F \to G) \to (1 \to G/\alpha(F))$ induces a canonical bijection $\mathbf{H}^0(F \to G) \xrightarrow{\sim} H^0(\operatorname{coker}\, \alpha)$.
\item If $\alpha\colon F \to G$ is surjective, then the embedding $(\ker \alpha \to 1) \hookrightarrow (F \to G)$ of complexes induces a canonical bijection $H^1(\ker \alpha) \xrightarrow{\sim} \mathbf{H}^0(F \to G)$.
\end{enumerate}

In the rest of this section we define the hypercohomology exact sequence assigned to a short exact sequence of complexes of $\Gamma$-groups.

\textbf{1.5} \textit{Exact sequences.} \quad A short exact sequence of complexes of $\Gamma$-groups is a sequence
\begin{equation*}
1 \to (F_1 \to G_1) \xrightarrow{i} (F_2 \to G_2) \xrightarrow{j} (F_3 \to G_3) \to 1
\end{equation*}
such that the rows in the commutative diagram
\begin{equation*}
\xymatrix{
1 \ar[r] & F_1 \ar[r] \ar[d]_{\alpha_1} & F_2 \ar[r] \ar[d]_{\alpha_2} & F_3 \ar[r] \ar[d]_{\alpha_3} & 1 \\
1 \ar[r] & G_1 \ar[r] & G_2 \ar[r] & G_3 \ar[r] & 1
}
\end{equation*}
are exact. We regard $F_1$ and $G_1$ as subgroups of $F_2$ and $G_2$, respectively. For such an exact sequence we define the connecting map
\begin{equation*}
\delta_{-1}\colon \mathbf{H}^{-1}(F_3 \to G_3) \to \mathbf{H}^0(F_1 \to G_1)
\end{equation*}
as follows.

Let $f_3 \in \mathbf{H}^{-1}(F_3 \to G_3) = (\ker \alpha_3)^\Gamma$. Choose $f \in F_2$ such that $f(\bmod\, F_1) = f_3$. We define a 0-cochain $(\varphi_1, g_1) \in C^0(F_2 \to G_2)$ by
\begin{equation*}
\varphi_1(\sigma) = f \cdot {}^{\upsig}f^{-1},\ g_1 = \alpha_2(f)
\end{equation*}
It is easy to show that $(\varphi_1, g_1) \in Z^0(F_1 \to G_1)$.

We set $\delta_{-1}(f_3) = \mathrm{Cl}(\varphi_1, g_1) \in \mathbf{H}^0(F_1 \to G_1)$. We leave to the reader to check that the map $\delta_{-1}$ is well-defined, i.e., $\delta_{-1}(f_3)$ does not depend on the choice of the representative $f \in F_2$ of $f_3$.

\subsection*{1.6 Proposition.} \textit{Let}
\begin{equation*}
1 \to (F_1 \to G_1) \xrightarrow{i} (F_2 \to G_2) \xrightarrow{j} (F_3 \to G_3) \to 1
\end{equation*}
\textit{be an exact sequence of complexes of $\Gamma$-groups. Then the hypercohomology sequence}
\begin{align}
&1 \longrightarrow \mathbf{H}^{-1}(F_1 \to G_1) \xrightarrow{i_*} \mathbf{H}^{-1}(F_2 \to G_2) \xrightarrow{j_*} \mathbf{H}^{-1}(F_3 \to G_3) \tag{1.6.1}\\
&\xrightarrow{\delta_{-1}} \mathbf{H}^0(F_1 \to G_1) \xrightarrow{i_*} \mathbf{H}^0(F_2 \to G_2) \xrightarrow{j_*} \mathbf{H}^0(F_3 \to G_3) \notag
\end{align}
\textit{is exact.}

Note that exactness makes sense because $\mathbf{H}^{-1}(F_k \to G_k)$ is a group and $\mathbf{H}^0(F_k \to G_k)$ is a pointed set $(k = 1, 2, 3)$.

\textit{Proof.} We prove the exactness at $\mathbf{H}^{-1}(F_3 \to G_3)$. It follows immediately from the definition of $\delta_{-1}$ that $\delta_{-1} \circ j_* = 1$. Conversely, suppose that $f_3 \in (\ker \alpha_3)^\Gamma$ and $\delta_{-1}(f_3) = 1$. Let $f$ be a representative of $f_3$ in $F_2$. Then there exists $f_1 \in F_1$ such that
\begin{align*}
\alpha_1(f_1)^{-1} \cdot \alpha_2(f) &= 1\\
f_1^{-1} \cdot f \cdot {}^{\upsig}\hm f^{-1} \cdot {}^{\upsig}\hm f_1 &= 1,
\end{align*}
hence
\begin{equation*}
\alpha_2(f_1^{-1}f) = 1,\ {}^{\upsig}(f_1^{-1}f) = f_1^{-1}f
\end{equation*}

Set $f' = f_1^{-1}f$. Then $f'(\mathrm{mod}\, F_1) = f_3$, $\alpha_2(f') = 1$, and ${}^{\upsig}\hm f' = f'$ for any $\sigma \in \Gamma$. Thus
\begin{equation*}
f' \in (\ker\, \alpha_2)^\Gamma = \mathbf{H}^{-1}(F_2 \to G_2)
\end{equation*}
and $f_3 = j_*(f')$. We have proved that $f_3 \in \mathrm{im}\, j_*$.

We leave the proof of the exactness at the other terms to the reader.
\qed

\section*{2 Crossed modules and $\mathbf{H}^1$}

To define $\mathbf{H}^1(F \to G)$ we need an additional structure on $F \to G$, namely the structure of \textit{crossed module}.

\textbf{2.1 Definition.} A crossed module is a short complex (homomorphism) $\alpha\colon F \to G$, endowed with a left action of $G$ on $F$ (denoted by $(g, f) \mapsto {}^\upg f$) satisfying
\begin{align}
ff'f^{-1} &= {}^{\alpha(f)\!}f' \tag{2.1.1}\\
\alpha({}^\upg f) &= g \cdot \alpha(f) \cdot g^{-1} \tag{2.1.2}
\end{align}
for any $f, f' \in F$, $g \in G$.

We say that a group $\Gamma$ acts on a crossed module $\alpha\colon F \to G$, if $\Gamma$ acts on $F$ and $G$ such that
\begin{equation*}
\alpha({}^{\upsig}\hm f) = {}^{\upsig}(\alpha(f)),\ \,
{}^{\upsig}({}^\upg f) = {}^{{}^\upsig g}({}^{\upsig}\hm f) \quad
\text{for any} \quad f \in F,\ g \in G,\ \sigma \in \Gamma.
\end{equation*}

\textbf{2.2 Examples} \textit{of crossed modules.}

\begin{enumerate}
\renewcommand{\labelenumi}{(\arabic{enumi})}
\item $\alpha\colon F \to G$ where $F$ is any (abelian) $G$-module, $\alpha$ is trivial.
\item $\alpha\colon F \hookrightarrow G$ where $F$ is a normal subgroup of $G$, $\alpha\colon F \hookrightarrow G$ is the inclusion, ${}^\upg f = gfg^{-1}$.
\item $\alpha\colon F \to G$ where $F \to G$ is any surjective homomorphism with central kernel. An element $g \in G$ acts on $F$ by ${}^\upg f = \tilde{g}f\tilde{g}^{-1}$ where $\tilde{g}$ is any lifting of $g$ to $F$.
\item $F \to \mathrm{Aut}\, F$ for any group $F$, $f \mapsto \mathrm{int}(f)$.
\item Let $X$ be a ``nice'' topological space, $Y \subset X$ a subspace and $x_0 \in Y$ a point. Then $\pi_1(Y, x_0)$ acts on $\pi_2(X, Y, x_0)$, and the complex $\pi_2(X, Y, x_0) \xrightarrow{\partial} \pi_1(Y, x_0)$ (where $\partial$ is the boundary homomorphism) is a crossed module.
\item Deligne's crossed module $\rho\colon G^{\mathrm{sc}} \to G$ of algebraic groups, described in the Introduction.
\end{enumerate}

\textbf{2.3} \textit{Remark.} J. H. C. Whitehead [W1], [W2], who introduced the notion of a crossed module, considered the crossed module 2.2(5). Dedecker showed in [Ded1], [Ded2] that
a crossed module $F\to G$ suits to define hypercohomology $\mathbf{H}^{1}(X,F\to G)$ where $X$ is a group, a topological space and so on. For a survey on crossed modules see [BHu].

\textbf{2.4. Lemma} (cf. [BHu]). \textit{Let $F\xrightarrow{\alpha}G$ be a crossed module. Then}
\begin{enumerate}
\renewcommand{\labelenumi}{(\roman{enumi})}
\item \textit{the group} $\ker\alpha$ \textit{is central in $F$;}
\item $\ker\alpha$ \textit{is $G$-invariant;}
\item $\operatorname{im}\alpha$ \textit{is normal in $G$.}
\end{enumerate}

\textit{Proof.} (i) follows from (2.1.1); (ii) and (iii) follow from (2.1.2).\qed

\textbf{2.5 Corollary.} \textit{The action of $G$ on $F$ induces an action of} $\operatorname{coker}\alpha$ \textit{on the abelian group} $\ker\alpha$.

\textbf{2.6} \textit{The group structure on $H^{0}$.} Let $F\to G$ be a crossed group with a $\Gamma$-action. We show that $C^{0}=C^{0}(F\to G)$, $Z^{0}(F\to G)$ and $\mathbf{H}^{0}(F\to G)$ have natural group structures.

The group $G$ acts on $\operatorname{Maps}(\Gamma,F)$ by $({}^{\upg}\varphi)(\sigma)={}^{\upg}(\varphi(\sigma))$ ($\varphi\in\operatorname{Maps}(\Gamma,F)$, $\sigma\in\Gamma$). We define a group structure on $C^{0}$ by
\begin{equation*}
(\varphi_{1},g_{1})\cdot(\varphi_{2},g_{2})=({\hs}^{g_{1}\!}\varphi_{2}\cdot\varphi_{1},\,g_{1}g_{2}).
\end{equation*}
One can check that $Z^{0}$ is a subgroup of $C^{0}$ with respect to this group structure.

Consider the map $\nu\colon F\to Z^{0}$ defined by the formula $\nu(f)=(\varphi,\alpha(f))$ where $\varphi(\sigma)=f\cdot{}^{\upsig}f^{-1}$. One can easily check that $\nu$ is a group homomorphism and its image is normal in $Z^{0}$. Moreover the right action of $F$ on $Z^{0}$ defined by
\begin{equation*}
((\varphi,g),f)\longmapsto\nu(f^{-1})\cdot(\varphi,g)
\end{equation*}
coincides with the action $*$ of 1.2. Thus $\mathbf{H}^{0}(F\to G)=Z^{0}/\operatorname{im}\nu$, and therefore $\mathbf{H}^{0}(F\to G)$ has a canonical group structure. This group structure depends functorially on the crossed module $F\to G$.

\textbf{2.7} \textit{Hypercohomology in degree $1$.} Let $F\to G$ be a crossed module with a $\Gamma$-action. Following Dedecker [Ded3] we define the first hypercohomology as follows.

Let $Z^{1}$ denote the set of pairs $(h,\psi)\in\operatorname{Maps}(\Gamma\times\Gamma,F)\times\operatorname{Maps}(\Gamma,G)$ such that for any $\sigma,\tau,\upsilon\in\Gamma$
\begin{align*}
\alpha(h(\sigma,\tau))^{-1}\cdot\psi(\sigma\tau) &= \cdot\psi(\sigma)\cdot{}^{\upsig}\psi(\tau)\\
h(\sigma,\tau\upsilon)\cdot{}^{\psi(\sigma)\sigma}h(\tau,\upsilon) &= h(\sigma\tau,\upsilon)\cdot h(\sigma,\tau).
\end{align*}
We define a right action $Z^{1}\times C^{0}\to Z^{1}$ of the group of $0$-cochains $C^{0}$ on the set of $1$-cocycles $Z^{1}$. For $(a,g)\in C^{0}$ we set
\begin{equation*}
(h,\psi)*(a,g)=(h',\psi')
\end{equation*}
where
\begin{align*}
\psi'(\sigma) &= g^{-1}\cdot\alpha(a(\sigma))\cdot\psi(\sigma)\cdot{}^{\upsig}g\\
h'(\sigma,\tau) &= {}^{g^{-1}}\!\left[a(\sigma\tau)\cdot h(\sigma,\tau)\cdot{\hs}^{\psi(\sigma)\sigma\!}a(\tau)^{-1}\cdot a(\sigma)^{-1}\right]
\end{align*}
One can easily check that this is a group action.

Now we set
\begin{equation*}
\mathbf{H}^{1}(F\to G)=Z^{1}/C^{0}.
\end{equation*}
The set $\mathbf{H}^{1}(F\to G)$ has \textit{a neutral element}, namely the class of the trivial cocycle $(1,1)\in Z^{1}$. We write $\mathrm{Cl}(h,\psi)$ for the hypercohomology class of $1$-cocycle $(h,\psi)$.

\textbf{2.8} \textit{Morphisms of crossed modules.} A morphism $\varepsilon\colon(F_{1}\to G_{1})\to(F_{2}\to G_{2})$ of crossed modules is a pair of homomorphisms $(\varepsilon_{0}\colon G_{1}\to G_{2}$, $\varepsilon_{-1}\colon F_{1}\to F_{2})$ such that the diagram
\begin{equation*}
\xymatrix{
F_1 \ar[r]^{\varepsilon_{-1}} \ar[d]_{\alpha_1} & F_2 \ar[d]^{\alpha_2} \\
G_1 \ar[r]_{\varepsilon_0} & G_2
}
\end{equation*}
commutes and ${}^{\varepsilon_{0}(g)}\varepsilon_{-1}(f)=\varepsilon_{-1}({}^{\upg}f)$ for any $g\in G_{1}$, $f\in F_{1}$.

A morphism $\varepsilon$ of crossed modules with $\Gamma$-action defines homomorphisms
\begin{equation*}
\varepsilon_{*}\colon\mathbf{H}^{i}(F_{1}\to G_{1})\to\mathbf{H}^{i}(F_{2}\to G_{2})\quad(i=-1,0)
\end{equation*}
and a map $\varepsilon_{*}\colon\mathbf{H}^{1}(F_{1}\to G_{1})\to\mathbf{H}^{1}(F_{2}\to G_{2})$ that takes the neutral element to the neutral element. Thus $\mathbf{H}^{-1}$, $\mathbf{H}^{0}$ and $\mathbf{H}^{1}$ are functors.

\textbf{2.9 Examples.}

\renewcommand{\labelenumi}{(\arabic{enumi})}
\begin{enumerate}
\item $\mathbf{H}^{1}(1\to G)=H^{1}(G)$.
\item $\mathbf{H}^{1}(F\to 1)=H^{2}(F)$ (note that in this case $F$ is abelian and therefore $H^{2}(F)$ makes sense). To $\mathrm{Cl}(h,1)\in\mathbf{H}^{1}(F\to 1)$
we assign $\mathrm{Cl}(h)\in H^{2}(F)$.
\item If $F\xrightarrow{\alpha}G$ is a crossed module and $\alpha$ is injective, then the morphism of complexes $(F\to G)\to(1\to G/\alpha(F))$ induces a canonical bijection $\mathbf{H}^{1}(F\to G)\overset{\sim}{\to}H^{1}(\operatorname{coker}\alpha)$.
\item If $\alpha$ is surjective, then the embedding $(\ker\alpha\to 1)\hookrightarrow(F\to G)$ of crossed modules induces a bijection $H^{2}(\ker\alpha)\overset{\sim}{\to}\mathbf{H}^{1}(F\to G)$. One can check that the map $H^{i}(G)\to\mathbf{H}^{i}(F\to G)=H^{i+1}(\ker\alpha)\ \ (i=0,1)$ coincides with the connecting map $\delta_{i}\colon H^{i}(G)\to H^{i+1}(\ker\alpha)$ assigned to the short exact sequence $1\to\ker\alpha\to F\to G\to 1$ (see [Se], Ch.\ I, \S5, for the definition of $\delta_{i}$).

In the rest of this section we prolong the hypercohomology exact sequence (1.6.1).
\end{enumerate}

\textbf{2.10} Let
\begin{equation*}
1\to(F_{1}\to G_{1})\stackrel{i}{\longrightarrow}(F_{2}\to G_{2})\stackrel{j}{\longrightarrow}(F_{3}\to G_{3})\to 1
\end{equation*}
be an exact sequence of complexes of groups. We identify $(F_{1}\to G_{1})$ with its image in $(F_{2}\to G_{2})$. Assume that $(F_{1}\to G_{1})$ and $(F_{2}\to G_{2})$ are endowed with structures of crossed modules such that $i$ is a morphism of crossed modules. We assume also that

(2.10.1) \quad $F_{1}$ \textit{is $G_{2}$-invariant in $F_{2}$.}

\noindent Then $G_{2}$ acts on $F_{3}\simeq F_{2}/F_{1}$. We do not assume that $(F_{3}\to G_{3})$ is a crossed module.

We define a left action of the group $\mathbf{H}^{0}(F_{2}\to G_{2})$ on the set $\mathbf{H}^{0}(F_{3}\to G_{3})$ by
\begin{equation*}
\mathrm{Cl}(\varphi_{2},g_{2})\cdot\mathrm{Cl}(\varphi_{3},g_{3})=\mathrm{Cl}({}^{g_{2}}\varphi_{3}\cdot j(\varphi_{2}),j(g_{2})\cdot g_{3})
\end{equation*}
One can check that this is a well-defined group action.

\textbf{2.11} \textit{The connecting map.} Let a short exact sequence
\begin{equation*}
1\to(F_{1}\to G_{1})\to(F_{2}\to G_{2})\to(F_{3}\to G_{3})\to 1
\end{equation*}
be as in 2.10. We define the connecting map
\begin{equation*}
\delta_{0}\colon\mathbf{H}^{0}(F_{3}\to G_{3})\to\mathbf{H}^{1}(F_{1}\to G_{1})
\end{equation*}
as follows.

Let $\xi_{3}\in\mathbf{H}^{0}(F_{3}\to G_{3})$, $\xi_{3}=\mathrm{Cl}(\varphi_{3},g_{3})$, $(\varphi_{3},g_{3})\in Z^{0}(F_{3}\to G_{3})$. We lift $(\varphi_{3},g_{3})$ to some $(\varphi,g)$, $\varphi\in\mathrm{Maps}(\Gamma,F_{2})$, $g\in G_{2}$. We set
\begin{align*}
\psi_{1}(\sigma) &= g^{-1}\cdot\alpha_{2}(\varphi(\sigma))\cdot{}^{\upsig}g\\
h_{1}(\sigma,\tau) &= {}^{g^{-1}}[\varphi(\sigma\tau)\cdot{}^{\upsig}\varphi(\tau)^{-1}\cdot\varphi(\sigma)^{-1}].
\end{align*}
Then $\psi_{1}(\sigma)\in G_{1}$ and $h_{1}(\sigma,\tau)\in F_{1}$ for any $\sigma,\tau\in\Gamma$ (we use (2.10.1)).

We show that $(h_{1},\psi_{1})\in Z^{1}(F_{1}\to G_{1})$. We have
\begin{align*}
\psi_{1}(\sigma)\cdot{}^{\upsig}\psi_{1}(\tau) &= g^{-1}\cdot\alpha_{2}(\varphi(\sigma))\cdot{}^{\upsig}g\cdot{}^{\upsig}g^{-1}\cdot\alpha_{2}({}^{\upsig}\varphi(\tau))\cdot{}^{\sigma\tau}g\\
&= g^{-1}\cdot\alpha_{2}(\varphi(\sigma)\cdot{}^{\upsig}\varphi(\tau)\cdot\varphi(\sigma\tau)^{-1})\cdot g\cdot g^{-1}\cdot\alpha_{2}(\varphi(\sigma\tau))\cdot{}^{\sigma\tau}g\\
&= \alpha_{2}(h_{1}(\sigma,\tau))^{-1}\cdot\psi_{1}(\sigma\tau);
\end{align*}
\begin{align*}
h_{1}(\sigma,\tau v)&\cdot{}^{\psi_{1}(\sigma)\sigma}h_{1}(\tau,v) \\
&= {}^{g^{-1}}[\varphi(\sigma\tau v)\cdot{}^{\upsig}\varphi(\tau v)^{-1}\cdot\varphi(\sigma)^{-1}]
     \cdot{}^{g^{-1}}[\varphi(\sigma)\cdot{}^{\upsig}\varphi(\tau v)\cdot{}^{\sigma\tau\!}\varphi(v)^{-1}\cdot{}^{\upsig}\varphi(\tau)^{-1}\cdot\varphi(\sigma)^{-1}]\\
&= {}^{g^{-1}}[\varphi(\sigma\tau v)\cdot{}^{\sigma\tau\!}\varphi(v)^{-1}\cdot{}^{\upsig}\varphi(\tau)^{-1}\cdot\varphi(\sigma)^{-1}]\\
&= {}^{g^{-1}}[\varphi(\sigma\tau v)\cdot{}^{\sigma\tau\!}\varphi(v)^{-1}\cdot\varphi(\sigma\tau)^{-1}]\cdot{}^{g^{-1}}[\varphi(\sigma\tau)\cdot{}^{\upsig}\varphi(\tau)^{-1}\cdot\varphi(\sigma)^{-1}]\\
&= h_{1}(\sigma\tau,v)\cdot h_{1}(\sigma,\tau).
\end{align*}
Hence $(h_{1},\psi_{1})\in Z^{1}(F_{1}\to G_{1})$.

We set $\delta_{0}(\xi_{3})=\mathrm{Cl}(h_{1},\psi_{1})\in\mathbf{H}^{1}(F_{1}\to G_{1})$. We leave to the reader to check that $\delta_{0}(\xi_{3})$ is well defined.

\subsection*{2.12 Proposition.} \textit{Let}
\begin{equation}
1 \rightarrow (F_1 \rightarrow G_1) \xrightarrow{i} (F_2 \rightarrow G_2) \xrightarrow{j} (F_3 \rightarrow G_3) \rightarrow 1 \tag{2.12.1}
\end{equation}
\textit{be a short exact sequence of complexes of} $\Gamma$\textit{-groups where} $i$ \textit{is an embedding of crossed modules with} $\Gamma$\textit{-action. We identify} $(F_1 \to G_1)$ \textit{with its image in} $(F_2 \to G_2)$ \textit{and assume that the subgroup} $F_1 \subset F_2$ \textit{is} $G_2$\textit{-invariant. Then}
\begin{enumerate}
\renewcommand{\labelenumi}{(\roman{enumi})}
\item \textit{the sequence}
\begin{align}
&1 \rightarrow \mathbf{H}^{-1}(F_1 \rightarrow G_1) \xrightarrow{i_*} \mathbf{H}^{-1}(F_2 \rightarrow G_2) \xrightarrow{j_*} \mathbf{H}^{-1}(F_3 \rightarrow G_3) \notag \\
\xrightarrow{\delta_{-1}}\ &\mathbf{H}^0(F_1 \rightarrow G_1) \xrightarrow{i_*} \mathbf{H}^0(F_2 \rightarrow G_2) \xrightarrow{j_*} \mathbf{H}^0(F_3 \rightarrow G_3) \tag{2.12.2} \\
\xrightarrow{\delta_0}\ &\mathbf{H}^1(F_1 \to G_1) \xrightarrow{i_*} \mathbf{H}^1(F_2 \to G_2) \notag
\end{align}
\textit{is exact.}
\item $\delta_0$ \textit{defines a bijection}
\begin{equation}
\mathbf{H}^0(F_2 \rightarrow G_2)\backslash\mathbf{H}^0(F_3 \rightarrow G_3) \xrightarrow{\sim} \ker[\mathbf{H}^1(F_1 \rightarrow G_1) \rightarrow \mathbf{H}^1(F_2 \rightarrow G_2)]. \tag{2.12.3}
\end{equation}
\end{enumerate}

\textit{Proof.} We leave the proof of (ii) to the reader. To prove (i), we must prove the exactness at the terms $\mathbf{H}^0(F_3 \to G_3)$ and $\mathbf{H}^1(F_1 \to G_1)$. We leave the proof for $\mathbf{H}^0(F_3 \to G_3)$ to the reader.

We prove the exactness at $\mathbf{H}^1(F_1 \to G_1)$. It is clear that $i_* \circ \delta_0 = 1$. Indeed, the cocycle $(h_1, \psi_1) \in Z^1(F_1 \to G_1)$ constructed in 2.11 is cohomologous to $(1,1)$ in $Z^1(F_2 \rightarrow G_2)$.

Conversely, let $\eta_1 \in \mathbf{H}^1(F_1 \to G_1)$, $\eta_1 = \operatorname{Cl}(h_1, \psi_1)$. Assume that $i_*(\eta_1) = 1$. Then
\begin{align*}
\psi_1(\sigma) &= g^{-1} \cdot \alpha_2(a(\sigma)) \cdot {}^{\upsig\hm} g \\
h_1(\sigma, \tau) &= {}^{g^{-1}}[a(\sigma\tau) \cdot {}^{\upsig} a(\tau)^{-1} \cdot a(\sigma)^{-1}]
\end{align*}
for some $a\colon \Gamma \to F_2$, $g \in G_2$. Set $g_3 = g(\mathrm{mod}\, G_1) \in G_3$, $\varphi_3(\sigma) = a(\sigma)(\mathrm{mod}\, F_1) \in F_3$. Using (2.10.1) one can easily check that $(\varphi_3, g_3) \in Z^0(F_3 \to G_3)$. Set $\xi_3 = \operatorname{Cl}(\varphi_3, g_3)$; then $\eta_1 = \delta_0(\xi_3)$. Thus $\eta_1 \in \operatorname{im}\, \delta_0$, which was to be proved. \qed

\medskip
\noindent\textbf{2.12.4} \textit{Remark.} The hypercohomology exact sequence (2.12.2) depends on the short exact sequence (2.12.1) functorially.

\medskip
\noindent\textbf{2.13 Corollary.} \textit{Let $F \xrightarrow{\alpha} G$ be a crossed module with $\Gamma$-action.}
\begin{enumerate}
\renewcommand{\labelenumi}{(\roman{enumi})}
\item \textit{{\rm ([Br1], (4.2.2))}. There is an exact sequence}
\begin{align}
1 \longrightarrow \mathbf{H}^{-1}(F \rightarrow G) \xrightarrow{\lambda_{-1}} H^0(F) \xrightarrow{\alpha_*} H^0(G) \xrightarrow{\varkappa_0} \mathbf{H}^0(F \rightarrow G) \tag{2.13.1} \\
\xrightarrow{\lambda_0} H^1(F) \xrightarrow{\alpha_*} H^1(G) \xrightarrow{\varkappa_1} \mathbf{H}^1(F \to G). \notag
\end{align}
\end{enumerate}

\noindent (ii) \textit{The map $\alpha_*\colon H^1(F) \to H^1(G)$ defines a bijection}
\begin{equation}
\mathbf{H}^0(F \to G)\backslash H^1(F) \xrightarrow{\sim} \ker[\varkappa_1\colon H^1(G) \to \mathbf{H}^1(F \to G)] \tag{2.13.2}
\end{equation}
Here $\alpha_*\colon H^i(F) \to H^i(G)$ are the canonical maps induced by $\alpha$. The maps $\lambda_{-1}$, $\varkappa_0$, $\lambda_0$ and $\varkappa_1$ can be described as follows:
\begin{align*}
\lambda_{-1}\colon \mathbf{H}^{-1}(F \to G) &= (\ker\,\alpha)^\Gamma \hookrightarrow F^\Gamma = H^0(F),\ f \longmapsto f \\
\varkappa_0\colon H^0(G) &= G^\Gamma \to \mathbf{H}^0(F \to G),\ g \longmapsto \operatorname{Cl}(1,g) \\
\lambda_0\colon \operatorname{Cl}(\varphi, g) &\longmapsto \operatorname{Cl}(\varphi) \\
\varkappa_1\colon \operatorname{Cl}(\psi) &\longmapsto \operatorname{Cl}(1,\psi)
\end{align*}

\medskip
\noindent\textbf{2.13.3} \textit{Remark.} The exact sequence (2.13.1), but without the last term, was earlier constructed by Deligne ([Del], (2.4.3.1)).

\medskip
\noindent\textit{Proof.} Consider the short exact sequence
\begin{equation*}
1 \to (1 \to G) \to (F \to G) \to (F \to 1) \to 1
\end{equation*}
of complexes of $\Gamma$-groups, where $(1 \to G) \to (F \to G)$ is a morphism of crossed modules. The exact sequence (2.12.2) takes in our case the form (2.13.1), and the bijection (2.12.3) takes the form (2.13.2). \qed

\medskip
\noindent\textbf{2.14} \textit{Twisting.} \quad To describe the fibers of the map $\varkappa_1\colon H^1(G) \to \mathbf{H}^1(F \to G)$ we need twisting.

The group $G$ acts on the crossed module $(F \to G)$. An element $g_* \in G$ acts by
\begin{equation*}
f \mapsto {}^{g_*\!}f,\ g \mapsto g_* g g_*^{-1}\ (f \in F,\, g \in G)
\end{equation*}

Let $\psi \in Z^1(G)$. We can define the twisted crossed module ${}_\psi(F \to G) = ({}_\psi F \to {}_\psi G)$, where the twisted groups ${}_\psi F$ and ${}_\psi G$ are the same $F$ and $G$ as abstract groups, but $\Gamma$ acts differently, namely,
\begin{equation*}
{}^{\sigma*\!}f = {}^{\psi(\sigma)\sigma\!\hm}f,\ {}^{\sigma*\!}g = \psi(\sigma) \cdot {}^{\upsig} g \cdot \psi(\sigma)^{-1}\ (\sigma \in \Gamma,\, f \in F,\, g \in G).
\end{equation*}

We define a map
\begin{equation*}
t_\psi\colon \mathbf{H}^1({}_\psi(F \to G)) \to \mathbf{H}^1(F \to G)
\end{equation*}
taking $1 \in \mathbf{H}^1({}_\psi(F \to G))$ to $\operatorname{Cl}(1,\psi) \in \mathbf{H}^1(F \to G)$. Let $(h', \psi') \in Z^1({}_\psi(F \to G))$. By definition this means that
\begin{align*}
\psi'(\sigma) \cdot \psi(\sigma) \cdot {}^{\upsig}\psi'(\tau) \cdot \psi(\sigma)^{-1} &= \alpha(h'(\sigma,\tau))^{-1} \cdot \psi'(\sigma\tau) \\
h'(\sigma, \tau v) \cdot {}^{\psi'(\sigma)\psi(\sigma)\sigma}h'(\tau, v) &= h'(\sigma\tau, v) \cdot h'(\sigma, \tau).
\end{align*}
We set
\begin{equation*}
t_\psi(\operatorname{Cl}(h', \psi')) = \operatorname{Cl}(h', \psi'\psi)
\end{equation*}
One can easily check that the map $t_\psi$ is well-defined.

We can define a map $t_\psi\colon H^1({}_\psi G) \to H^1(G)$ in a similar way. The diagram
\begin{equation}\tag{2.14.1}
\begin{aligned}
\xymatrix@C=15mm{
H^1({}_\psi G) \ar[r]^{t_\psi} \ar[d]_-{{}_\psi\varkappa_1} & H^1(G) \ar[d]^{\varkappa_1} \\
\mathbf{H}^1({}_\psi(F \to G)) \ar[r]^-{t_\psi} & \mathbf{H}^1(F \to G)
}
\end{aligned}
\end{equation}
commutes.

\medskip
\noindent\textbf{2.15 Proposition.} \textit{Let $(F \xrightarrow{\alpha} G)$ be a crossed module with $\Gamma$-action. Consider the exact sequence (2.13.1). Let $\eta \in H^1(G)$, $\eta = \operatorname{Cl}(\psi)$, $\psi \in Z^1(G)$. Then the fiber of $\varkappa_1$ over $\varkappa_1(\eta)$ is in canonical bijection with the quotient set}
\begin{equation*}
\mathbf{H}^0({}_\psi(F \to G))\backslash H^1({}_\psi F).
\end{equation*}

\textit{Proof.} The map $t_\psi\colon \mathbf{H}^1({}_\psi(F \to G)) \to \mathbf{H}^1(F \to G)$ takes $1$ to $\operatorname{Cl}(1,\psi) = \varkappa_1(\eta)$. Since the diagram (2.14.1) is commutative, the map $t_\psi\colon H^1({}_\psi G) \to H^1(G)$ takes the kernel of ${}_\psi\varkappa_1$ to the fiber of $\varkappa_1$ over $\varkappa_1(\eta)$. By Corollary 2.13 (ii) the kernel of ${}_\psi\varkappa_1$ is in canonical bijection with $\mathbf{H}^0({}_\psi(F \to G))\backslash H^1({}_\psi F)$. This proves the proposition. \qed

\medskip
\noindent\textbf{2.16 Proposition} ([Br1], (5.1.3)). \textit{Let}
\begin{equation*}
1 \to (F_1 \to G_1) \xrightarrow{i} (F_2 \to G_2) \xrightarrow{j} (F_3 \to G_3) \to 1
\end{equation*}
\textit{be an exact sequence of crossed modules with $\Gamma$-action. Then the sequence}
\begin{align}
&1 \longrightarrow \mathbf{H}^{-1}(F_1 \to G_1) \xrightarrow{i_*} \mathbf{H}^{-1}(F_2 \to G_2) \xrightarrow{j_*} \mathbf{H}^{-1}(F_3 \to G_3) \notag \\
\xrightarrow{\delta_{-1}}\ &\mathbf{H}^0(F_1 \to G_1) \xrightarrow{i_*} \mathbf{H}^0(F_2 \to G_2) \xrightarrow{j_*} \mathbf{H}^0(F_3 \to G_3) \tag{2.16.1} \\
\xrightarrow{\delta_0}\ &\mathbf{H}^1(F_1 \to G_1) \xrightarrow{i_*} \mathbf{H}^1(F_2 \to G_2) \xrightarrow{j_*} \mathbf{H}^1(F_3 \to G_3) \notag
\end{align}
\textit{is defined and exact.}

\textit{Proof.} Since $j$ is a morphism of crossed modules, the subgroup $F_1 \subset F_2$ is $G_2$-invariant, and therefore the map $\delta_0$ is defined. We must prove only the exactness at $\mathbf{H}^1(F_2 \to G_2)$; we leave it to the reader. \qed

\medskip
\noindent\textbf{2.17} \textit{The case of a normal abelian submodule.} We wish to extend the exact sequence (2.16.1). We assume that the crossed submodule $(F_1 \to G_1) \subset (F_2 \to G_2)$ is \textit{abelian}, i.e., $F_1$ and $G_1$ are abelian groups and $G_1$ acts on $F_1$ trivially. We assume also that

(2.17.1) $\alpha_{2}(F_{2})$ \textit{commutes with} $G_{1}$ \textit{in} $G_{2}$;

(2.17.2) $F_{1}$ \textit{is central in} $F_{2}$.

It follows from (2.17.1) and (2.17.2) that the group $G_{3}$ acts on the complex $(F_{1}\to G_{1})$ through $\operatorname{coker}\alpha_{3}$. A cocycle $(h_{3},\psi_{3})\in Z^{1}(F_{3}\to G_{3})$ defines a cocycle $\bar{\psi}_{3}\in Z^{1}(\operatorname{coker}\alpha_{3})$, namely $\bar{\psi}_{3}(\sigma)=\psi_{3}(\sigma)(\operatorname{mod}\alpha_{3}(F_{3}))$. Since $\operatorname{coker}\alpha_{3}$ acts on the complex $(F_{1}\to G_{1})$, we can define the twisted complex $_{\bar{\psi}_{3}}(F_{1}\to G_{1})$. We write $_{\psi_{3}}(F_{1}\to G_{1})$ for $_{\bar{\psi}_{3}}(F_{1}\to G_{1})$.

We define a hypercohomology class $\Delta_{1}(h_{3},\psi_{3})\in\mathbf{H}^{2}({}_{\psi_{3}}(F_{1}\to G_{1}))$ as follows. We lift $\psi_{3}$ to some continuous map $\psi\colon\Gamma\to G_{2}$ and lift $h_{3}$ to some continuous map $h\colon\Gamma\times\Gamma\to F_{2}$. Then we set
\begin{align*}
d_{1}(\sigma,\tau) &= \psi(\sigma)\cdot {}^{\upsig}\psi(\tau)\cdot\psi(\sigma\tau)^{-1}\cdot\alpha_{2}(h(\sigma,\tau))\\
 {\mathfrak{a}_{1}}(\sigma,\tau,\upsilon) &= {}^{\psi(\sigma)\sigma}h(\tau,\upsilon)^{-1}\cdot h(\sigma,\tau\upsilon)^{-1}\cdot h(\sigma\tau,\upsilon)\cdot h(\sigma,\tau)
\end{align*}
It is clear that $d_{1}(\sigma,\tau)\in G_{1}$, $ {\mathfrak{a}_{1}}(\sigma,\tau,\upsilon)\in F_{1}$. We must show now that $( {\mathfrak{a}_{1}},d_{1})\in Z^{2}({}_{\psi_{3}}(F_{1}\to G_{1}))$, i.e.,
\begin{gather*}
{}^{\psi(\sigma)\sigma}d_{1}(\tau,\upsilon)^{-1}d_{1}(\sigma,\tau)d_{1}(\sigma\tau,\upsilon)d_{1}(\sigma,\tau\upsilon)^{-1} = \alpha_{1}( {\mathfrak{a}_{1}}(\sigma,\tau,\upsilon)),\\
{}^{\psi(\sigma)\sigma} {\mathfrak{a}_{1}}(\tau,\upsilon,\rho)\cdot {\mathfrak{a}_{1}}(\sigma\tau,\upsilon,\rho)^{-1}\cdot {\mathfrak{a}_{1}}(\sigma,\tau\upsilon,\rho)\cdot {\mathfrak{a}_{1}}(\sigma,\tau,\upsilon\rho)^{-1} {\mathfrak{a}_{1}}(\sigma,\tau,\upsilon) = 1
\end{gather*}
We skip this tedious (though not easy) calculation.

We set $\Delta_{1}(h_{3},\psi_{3})=\operatorname{Cl}( {\mathfrak{a}_{1}},d_{1})\in\mathbf{H}^{2}({}_{\psi_{3}}(F_{1}\to G_{1}))$. We must check that the cohomology class $\delta_{1}(h_{3},\psi_{3})$ is well-defined, i.e., it does not depend on the choice of the lifting $(h,\psi)$ of $(h_{3},\psi_{3})$. We leave the check to the reader.

\textbf{2.18 Proposition.} \textit{Let}
\begin{equation*}
1\to(F_{1}\to G_{1})\overset{i}{\longrightarrow}(F_{2}\to G_{2})\overset{j}{\longrightarrow}(F_{3}\to G_{3})\to 1
\end{equation*}
\textit{be an exact sequence of crossed modules with $\Gamma$-action such that the crossed module $(F_{1}\to G_{1})$ is abelian and $(2.17.1)$ and $(2.17.2)$ hold. Let $(h_{3},\psi_{3})\in Z^{1}(F_{3}\to G_{3})$. Then $\operatorname{Cl}(h_{3},\psi_{3})\in\operatorname{im}j_{*}$ if and only if $\Delta_{1}(h_{3},\psi_{3})=1$.}

\textit{Proof.} Left to the reader. \qed

\textbf{2.19} \textit{The fibers of $j_{*}$.} Let the exact sequence
\begin{equation*}
1\to(F_{1}\to G_{1})\overset{i}{\longrightarrow}(F_{2}\to G_{2})\overset{j}{\longrightarrow}(F_{3}\to G_{3})\to 1
\end{equation*}
be as in Proposition 2.18. We want to describe the fibers of the map $j_{*}\colon\mathbf{H}^{1}(F_{2}\to G_{2})\to\mathbf{H}^{1}(F_{3}\to G_{3})$.

Let
\begin{align*}
(h,\psi)&\in Z^{1}(F_{2}\to G_{2}),\ \ \eta_{2}=\mathrm{Cl}(h,\psi)\in\mathbf{H}^{1}(F_{2}\to G_{2}),\\
(h_{3},\psi_{3})&=j(h,\psi),\ \ \eta_{3}=j_{*}(\eta_{2})=\mathrm{Cl}(h_{3},\psi_{3})\in\mathbf{H}^{1}(F_{3}\to G_{3}).
\end{align*}
We define a map $t_{(h,\psi)}\colon\mathbf{H}^{1}({}_{\psi_{3}}(F_{1}\to G_{1}))\to\mathbf{H}^{1}(F_{2}\to G_{2})$ which takes $1$ to $\eta_{2}$. We set
\begin{equation*}
t_{(h,\psi)}(\mathrm{Cl}(h_{1},\psi_{1}))=\mathrm{Cl}(hh_{1},\psi_{1}\psi)
\end{equation*}
One can check that $(hh_{1},\psi_{1}\psi)\in Z^{1}(F_{2}\to G_{2})$ and that the map $t_{(h,\psi)}$ is  well-defined.

\textbf{2.19.1 Lemma.} \textit{The fiber of the map $j_{*}\colon\mathbf{H}^{1}(F_{2}\to G_{2})\to\mathbf{H}^{1}(F_{3}\to G_{3})$ over $\eta_{3}=j_{*}(\eta_{2})$ is the image of the map $t_{(h,\psi)}$.}

\textit{Proof.} Easy. \qed

\textbf{2.20 Example.} Let $F\xrightarrow{\alpha}G$ be a crossed module with $\Gamma$-action. Consider the canonical exact sequence of crossed modules
\begin{equation*}
\text{(2.20.1)}\qquad\qquad 1\to(\ker\alpha\to 1)\xrightarrow{i}(F\to G)\xrightarrow{j}(F/\ker\alpha\hookrightarrow G)\to 1
\end{equation*}
The complex $(\ker\alpha\to 1)$ is abelian, and
\begin{align*}
\mathbf{H}^{i}(\ker\alpha\to 1)&=H^{i+1}(\ker\alpha)\quad(i\geq-1),\\
\mathbf{H}^{i}(F/\ker\alpha\hookrightarrow G)&=H^{i}(\operatorname{coker}\alpha)\quad\text{for }i=0,1.
\end{align*}
Note that conditions (2.17.1) and (2.17.2) are satisfied. By Propositions 2.16 and 2.18, to the short exact sequence (2.20.1) we can assign the hypercohomology exact sequence
\begin{align*}
\text{(2.20.2)}\qquad\quad 1\longrightarrow H^{1}(\ker\alpha)&\xrightarrow{i_{*}}\mathbf{H}^{0}(F\to G)\xrightarrow{j_{*}}(\operatorname{coker}\alpha)^{\Gamma}\xrightarrow{\delta_{0}}H^{2}(\ker\alpha)\\
&\longrightarrow\mathbf{H}^{1}(F\to G)\xrightarrow{j_{*}}H^{1}(\operatorname{coker}\alpha)\text{- - - - $\succ H^{3}({}_{\psi_{3}}(\ker\alpha))$}.
\end{align*}

The arrow $\text{- - - -}\succ$ in (2.20.2) is not a map, it just indicates that if $\eta_{3}\in H^{1}(\operatorname{coker}\alpha)$, $\eta_{3}=\mathrm{Cl}(\psi_{3})$, where $\psi_{3}\in Z^{1}(\operatorname{coker}\alpha)$, then $\eta_{3}$ comes from $\mathbf{H}^{1}(F\to G)$ if and only if $\Delta_{1}(\psi_{3})=1$. The fiber $j_{*}^{-1}(\eta_{3})$ is described in Lemma 2.19.1

\textbf{2.21} \textit{The case of a central submodule.} Let
\begin{equation*}
1\to(F_{1}\to G_{1})\xrightarrow{i}(F_{2}\to G_{2})\xrightarrow{j}(F_{3}\to G_{3})\to 1
\end{equation*}
be a short exact sequence. We identify the crossed module $(F_{1}\to G_{1})$ with its image in $(F_{2}\to G_{2})$.

We say that $(F_1 \to G_1)$ is \textit{central} in $(F_2 \to G_2)$, if $G_1$ is central in $G_2$, $F_1$ is central
in $F_2$, and $G_2$ acts trivially on $F_1$. Assume that $(F_1 \to G_1)$ is central in $(F_2 \to G_2)$. Then
we can define the connecting map $\delta_1\colon \mathbf{H}^1(F_3 \to G_3) \to \mathbf{H}^2(F_1 \to G_1)$.

Let $\eta_3 \in \mathbf{H}^1(F_3 \to G_3)$, $\eta_3 = \operatorname{Cl}(h_3, \psi_3)$. Then $\Delta_1(h_3, \psi_3) \in \mathbf{H}^2(F_1 \to G_1)$ (we write
$\mathbf{H}^2(F_1 \to G_1)$ instead of $\mathbf{H}^2({}_{\psi_3}(F_1 \to G_1))$ because $(F_1 \to G_1)$ is central in $(F_2 \to G_2)$\,).
One can check that $\Delta_1(h_3, \psi_3)$ does not depend on the choice of the cocycle $(h_3, \psi_3)$
representing $\eta_3$. We set
\begin{equation*}
\delta_1(\eta_3) = \Delta_1(h_3, \psi_3).
\end{equation*}

Propositions 2.16 and 2.18 imply

\noindent\textbf{2.22 Proposition.} \textit{Let}
\begin{equation*}
1 \to (F_1 \to G_1) \xrightarrow{i} (F_2 \to G_2) \xrightarrow{j} (F_3 \to G_3) \to 1
\end{equation*}
\textit{be a short exact sequence of crossed modules with $\Gamma$-action, where the crossed submodule
$(F_1 \to G_1)$ is central in $(F_2 \to G_2)$. Then the sequence}
\noindent(2.22.1) \hspace{1em} $\mathbf{H}^1(F_1 \to G_1) \xrightarrow{i_*} \mathbf{H}^1(F_2 \to G_2) \xrightarrow{j_*} \mathbf{H}^1(F_3 \to G_3) \xrightarrow{\delta_1} \mathbf{H}^2(F_1 \to G_1)$
\noindent\textit{is exact.}

\section*{3 Quasi-isomorphisms}

Let $(F_1 \xrightarrow{\alpha_1} G_1) \to (F_2 \xrightarrow{\alpha_2} G_2)$ be a morphism of crossed modules. Such a morphism
induces group homomorphisms $\ker \alpha_1 \longrightarrow \ker \alpha_2$, $\operatorname{coker} \alpha_1 \longrightarrow \operatorname{coker} \alpha_2$.

\medskip

\noindent\textbf{3.1 Definition.} A morphism $(F_1 \xrightarrow{\alpha_1} G_1) \to (F_2 \xrightarrow{\alpha_2} G_2)$ is called \textit{a quasi-isomorphism} if
the induced homomorphisms $\ker \alpha_1 \to \ker \alpha_2$ and $\operatorname{coker} \alpha_1 \to \operatorname{coker} \alpha_2$ are isomorphisms.

\medskip

\noindent\textbf{3.2 Examples.}

\renewcommand{\labelenumi}{(\arabic{enumi})}
\begin{enumerate}
\item Let $(F \xrightarrow{\alpha} G)$ be a crossed module. If $\alpha$ is injective, then $(F \to G) \to (1 \to \operatorname{coker} \alpha)$
is a quasi-isomorphism. If $\alpha$ is surjective, then $(\ker \alpha \to 1) \to (F \to G)$ is a quasi-isomorphism.

\item The morphism of crossed modules of algebraic groups $(Z^{(\mathrm{sc})} \to Z) \hookrightarrow (G^{\mathrm{sc}} \to G)$,
described in the introduction, is a quasi-isomorphism.
\end{enumerate}

\noindent\textbf{3.3 Theorem.} \textit{Let $\varepsilon\colon (F_1 \to G_1) \to (F_2 \to G_2)$ be a quasi-isomorphism of crossed modules
with $\Gamma$-action. Then $\varepsilon$ induces bijections}
\begin{equation*}
\varepsilon_*\colon \mathbf{H}^i(F_1 \to G_1) \longrightarrow \mathbf{H}^i(F_2 \to G_2)
\end{equation*}
\textit{for $i = -1, 0, 1$.}

\noindent\textit{Proof.} For $i = -1$ the assertion is obvious.

Let $i = 0$. From 2.20 we obtain a commutative diagram
\begin{equation*}
\xymatrix@C=1.45em{
1 \ar[r] & H^1(\ker \alpha_1) \ar[r] \ar[d]_{\sim} & \mathbf{H}^0(F_1 \to G_1) \ar[r] \ar[d]^{\varepsilon_*} & (\operatorname{coker} \alpha_1)^\Gamma \ar[r] \ar[d]_{\sim} & H^2(\ker \alpha_1) \ar[d]_{\sim} \\
1 \ar[r] & H^1(\ker \alpha_2) \ar[r] & \mathbf{H}^0(F_2 \to G_2) \ar[r] & (\operatorname{coker} \alpha_2)^\Gamma \ar[r] & H^2(\ker \alpha_2)
}
\end{equation*}
with exact rows. Three vertical arrows in this diagram are isomorphisms because $\varepsilon$ is a
quasi-isomorphism. Then by the five-lemma the map $\varepsilon_*\colon \mathbf{H}^0(F_1 \to G_1) \to \mathbf{H}^0(F_2 \to G_2)$
is also an isomorphism, which was to be proved.

Let $i = 1$. From 2.20 we obtain a commutative diagram
\begin{equation*}
\xymatrix@C=0.65em{
(\operatorname{coker} \alpha_1)^\Gamma \ar[r] \ar[d]_{\sim} & H^2(\ker \alpha_1) \ar[r] \ar[d]_{\sim} & \mathbf{H}^1(F_1 \to G_1) \ar[r] \ar[d]^{\varepsilon_*} & H^1(\operatorname{coker} \alpha_1) \ar@{-->}[r] \ar[d]_{\sim} & H^3({}_{\psi^{(1)}_3}(\ker \alpha_1)) \ar[d]_{\sim} \\
(\operatorname{coker} \alpha_2)^\Gamma \ar[r] & H^2(\ker \alpha_2) \ar[r] & \mathbf{H}^1(F_2 \to G_2) \ar[r] & H^1(\operatorname{coker} \alpha_2) \ar@{-->}[r] & H^3({}_{\psi^{(2)}_3}(\ker \alpha_2))
}
\end{equation*}
with exact rows. Four vertical arrows in this diagram are bijections because $\varepsilon$ is a quasi-isomorphism. We prove the assertion by diagram chasing. To prove the surjectivity of the
map $\varepsilon_*\colon \mathbf{H}^1(F_1 \to G_1) \to \mathbf{H}^1(F_2 \to G_2)$ we use Lemma 2.19.1. \qed

\section*{4 Abelianization maps}

Let $K$ be a field of characteristic 0, and $\bar{K}$ an algebraic closure of $K$. We set
$\Gamma = \operatorname{Gal}(\bar{K}/K)$.

The notions of a crossed module of algebraic groups and a quasi-isomorphism of
crossed modules of algebraic groups are defined in the evident way. If $F \to G$ is a crossed
module of algebraic groups, then $F(\bar{K}) \to G(\bar{K})$ is a (discrete) crossed module with a
$\Gamma$-action. We define the Galois hypercohomology of $F \to G$ by
\begin{equation*}
\mathbf{H}^i(K, F \to G) = \mathbf{H}^i(\Gamma, F(\bar{K}) \to G(\bar{K})) \qquad (i = -1, 0, 1).
\end{equation*}
We often abbreviate $\mathbf{H}^i(K, F \to G)$ to $\mathbf{H}^i(F \to G)$.

If $(F_1 \to G_1) \to (F_2 \to G_2)$ is a quasi-isomorphism of crossed modules of $K$-groups,
then
\begin{equation*}
\left(F_1(\bar{K}) \to G_1(\bar{K})\right) \to \left(F_2(\bar{K}) \to G_2(\bar{K})\right)
\end{equation*}
is a quasi-isomorphism of crossed modules with $\Gamma$-action, and by Theorem 3.3 we have a
bijection $\mathbf{H}^i(F_1 \to G_1) \xrightarrow{\sim} \mathbf{H}^i(F_2 \to G_2)$.

\noindent\textbf{4.1} Let $G$ be a connected reductive $K$-group. Let $G^{\mathrm{ss}}$ denote its derived group (which is
semisimple), and let $G^{\mathrm{sc}} \to G^{\mathrm{ss}}$ be the universal covering of $G^{\mathrm{ss}}$. Consider the composition
\begin{equation*}
\rho\colon G^{\mathrm{sc}} \to G^{\mathrm{ss}} \to G.
\end{equation*}
Then $G$ acts on $G^{\mathrm{sc}}$, and $G^{\mathrm{sc}} \overset{\rho}{\longrightarrow} G$ is a crossed module of $K$-groups. Let $Z$ denote the
center of $G$, and $Z^{(\mathrm{sc})}$ the center of $G^{\mathrm{sc}}$.

Let $T \subset G$ be a maximal torus defined over $K$. We set $T^{(\mathrm{sc})} = \rho^{-1}(T)$. We define the
abelian Galois cohomology $H^i_{\mathrm{ab}}(K, G)$ (which we usually abbreviate to $H^i_{\mathrm{ab}}(G)$) by
\begin{equation*}
H^i_{\mathrm{ab}}(K, G) := \mathbf{H}^i(K, T^{(\mathrm{sc})} \to T) = \mathbf{H}^i(K, Z^{(\mathrm{sc})} \to Z) \qquad (i \geq -1),
\end{equation*}
where we identify the abelian groups $\mathbf{H}^i(K, T^{(\mathrm{sc})} \to T)$ and $\mathbf{H}^i(K, Z^{(\mathrm{sc})} \to Z)$ using the
quasi-isomorphism $(Z^{(\mathrm{sc})} \to Z) \to (T^{(\mathrm{sc})} \to T)$ of abelian complexes. Note that $\mathbf{H}^i_{\mathrm{ab}}(K, \cdot)$
is a functor from the category of connected reductive $K$-group to the category of abelian
groups. We are interested here in $H^0_{\mathrm{ab}}$ and $H^1_{\mathrm{ab}}$.

\medskip

\noindent\textbf{4.1.1 Lemma.} \textit{For $i = 0, 1$ there is a canonical and functorial in $G$ bijection $\mathbf{H}^i(G^{\mathrm{sc}} \to
G) \overset{\sim}{\longrightarrow} H^i_{\mathrm{ab}}(G)$, which is a group isomorphism when $i = 0$.}

\noindent\textit{Proof.} The assertion follows from Theorem 3.3 applied to the quasi-isomorphisms
\begin{equation*}
(Z^{(\mathrm{sc})} \to Z) \to (T^{(\mathrm{sc})} \to T) \to (G^{\mathrm{sc}} \to G).\qquad\qed
\end{equation*}

\noindent\textbf{4.1.2} \textit{Remark} (essentially due to L. Breen). There is another, more intrinsic explanation
of the fact that $\mathbf{H}^1(G^{\mathrm{sc}} \to G)$ has a canonical structure of abelian group. Deligne ([De],
2.0.2) noted that the commutator morphism
\begin{equation*}
(g_1, g_2) \mapsto g_1 g_2 g_1^{-1} g_2^{-1} : \qquad G \times G \to G
\end{equation*}
can be uniquely lifted to a morphism
\begin{equation*}
(g_1, g_2) \mapsto \{g_1, g_2\} : \quad G \times G \to G^{\mathrm{sc}},
\end{equation*}
and we have $\{g_1, g_2\} = \{g_2, g_1\}$. The crossed module $G^{\mathrm{sc}} \to G$ together with the map
$\{\,,\,\}$ is a stable crossed module in the terminology of Conduch\'{e} ([Co], 3.1). To the crossed
module $G^{\mathrm{sc}} \to G$ one assigns a (fibered) Picard category $\mathcal{C}(G^{\mathrm{sc}} \to G)$ (cf. [Br2], the
remark after Def. 1.1.6). The map $\{\,,\,\}$ defines a commutativity constraint in $\mathcal{C}(G^{\mathrm{sc}} \to G)$,
and thus turns it into a commutative Picard category (cf. [Br3]). A commutative Picard
category is a categoric analogue of an abelian group, so the set of isomorphism classes of
torsors under such a category has a canonical structure of abelian group. Since $\mathbf{H}^1(G^{\mathrm{sc}} \to
G)$ is the set of isomorphism classes of torsors under $\mathcal{C}(G^{\mathrm{sc}} \to G)$ ([Br1], 6.2), we see that
$\mathbf{H}^1(G^{\mathrm{sc}} \to G)$ has a canonical structure of abelian group.

\medskip

\noindent\textbf{4.2} For $i = 0, 1$, we define the abelianization map $\mathrm{ab}^i$ as the composition
\begin{equation*}
\mathrm{ab}^i\colon H^i(G) \overset{\varkappa_i}{\longrightarrow} \mathbf{H}^i(G^{\mathrm{sc}} \to G) 
   \overset{\sim}{\longrightarrow} \mathbf{H}^i(Z^{(\mathrm{sc})} \to Z) = H^i_{\mathrm{ab}}(G),
\end{equation*}
where the map $\varkappa_i$ is induced by the embedding $(1 \to G) \to (G^{\mathrm{sc}} \to G)$ of crossed modules.
By Corollary 2.13(i) we have an exact sequence

\noindent(4.2.1) \hspace{0.5em} $G^{\mathrm{sc}}(K) \overset{\rho_*}{\longrightarrow} G(K) \overset{\mathrm{ab}^0}{\longrightarrow} H^0_{\mathrm{ab}}(K,G) \longrightarrow H^1(K, G^{\mathrm{sc}}) \overset{\rho_*}{\longrightarrow} H^1(K,G) \overset{\mathrm{ab}^1}{\longrightarrow} H^1_{\mathrm{ab}}(K,G)$

In [Bo3] (see also [Bo1]) we prove:
\subsection*{4.2.2 Proposition.} \textit{If $K$ is a local field of characteristic $0$ (archimedean or not) or a number field, then the map $\mathrm{ab}^{1}$ is surjective.}

\medskip
\noindent From Proposition 4.2.2 we deduce here

\subsection*{4.2.3 Corollary.} \textit{If $K$ is a non-archimedean local field of characteristic $0$, then the map $\mathrm{ab}^{1}$ is bijective.}

\noindent \textit{Proof.} Consider the exact sequence (4.2.1). By Proposition 2.15 any fiber of the map $\mathrm{ab}^{1}$ comes from $H^{1}(K, {}_{\psi}G^{\mathrm{sc}})$ where $\psi \in Z^{1}(K,G)$. Since ${}_{\psi}G^{\mathrm{sc}}$ is simply connected, by Kneser's theorem ([Kn]) we have $H^{1}(K, {}_{\psi}G^{\mathrm{sc}}) = 1$, so the map $\mathrm{ab}^{1}$ is injective. By Proposition 4.2.2 the map $\mathrm{ab}^{1}$ is surjective. We conclude that the abelianization map $\mathrm{ab}^{1}$ is bijective, which was to be proved.\qed

\noindent We see that when $K$ is a non-archimedean local field, the set $H^{1}(K,G)$ has a canonical and functorial structure of abelian group. (This result is due to Kottwitz [Ko1], [Ko2] in a slightly less functorial form.)

\medskip
\noindent \textbf{4.3} We can now describe the abelianization maps
\begin{equation*}
\mathrm{ab}^{i}\colon H^{i}(K,G) \to H^{i}_{\mathrm{ab}}(K,G) = \mathbf{H}^{i}(K, Z^{(\mathrm{sc})} \to Z) \qquad (i = 0,1)
\end{equation*}
\noindent explicitly in terms of cocycles.

\medskip
\noindent \textbf{4.3.1 Proposition.} \textit{Let $g \in H^{0}(k,G) = G(K)$. Write $g = \rho(g') \cdot z$ where $g' \in G^{\mathrm{sc}}(\bar{K})$, $z \in Z(\bar{K})$. Then $\mathrm{ab}^{0}(g) = \mathrm{Cl}(\varphi, z)$ where the map $\varphi\colon \Gamma \to Z^{(\mathrm{sc})}(\bar{K})$ is defined by $\varphi(\sigma) = (g')^{-1} \cdot {}^{\upsig}g'$.}

\noindent \textit{Proof.} We have $(1,g) * g' = (\varphi, z)$ with the notation of 1.2. Thus the $0$-cocycles $(1,g)$ and $(\varphi, z)$ are cohomological in $\mathbf{H}^{0}(K, G^{\mathrm{sc}} \to G)$. This proves the assertion.\qed

\medskip
\noindent \textbf{4.3.2 Proposition.} \textit{Let $\xi \in H^{1}(K,G)$ be a cohomology class, $\xi = \mathrm{Cl}(\psi)$, $\psi \in Z^{1}(K,G)$. Write $\psi(\sigma) = \rho(\psi'(\sigma)) \cdot z(\sigma)$ for $\sigma \in \Gamma$, where $\psi'\colon \Gamma \to G^{\mathrm{sc}}(\bar{K})$ and $z\colon \Gamma \to Z(\bar{K})$ are continuous maps. Then $\mathrm{ab}^{1}(\xi) = \mathrm{Cl}(h,z)$, where the map $h\colon \Gamma \times \Gamma \to Z^{(\mathrm{sc})}(\bar{K})$ is given by}
\begin{equation*}
h(\sigma, \tau) = \psi'(\sigma) \cdot {}^{\upsig}\psi'(\tau) \cdot \psi'(\sigma\tau)^{-1}.
\end{equation*}

\noindent \textit{Proof.} With the notation of 2.7 we have $(1,\psi) * ((\psi')^{-1}, 1) = (h,z)$. Thus the $1$-cocycles $(1,\psi)$ and $(h,z)$ are cohomological in $\mathbf{H}^{1}(K, G^{\mathrm{sc}} \to G)$. This proves the assertion.\qed

\medskip
\noindent {\sc Acknowledgements.} This text emerged as an attempt to understand cohomological constructions of Kottwitz [Ko2].
I found the necessary techniques in Dedecker [Ded3].

\noindent \hspace*{1.5em} The text was conceived during my stay at the Institute of Information Transmission of the USSR Academy of Sciences (Moscow),
written at the Institute for Advanced Study (Princeton),
and revised at the Max-Planck-Institut f\"{u}r Mathematik (Bonn). I am grateful to these institutes for their hospitality and support.

\noindent \hspace*{1.5em} I am very grateful to L. Breen for a number of useful remarks on the first version of this text.
It is a pleasure to thank P. Deligne for numerous helpful discussions.
Last but not least, I am deeply grateful to R.E. Kottwitz who explained to me in the summer of 1989 that the group I was interested in, was an (abelian) {\em hypercohomology} group.

\medskip
\noindent {\sc Disclosure.}
This text with formulas was converted from PDF into LaTeX using AI Underleaf.
The commutative diagrams were converted  into XY-pic using ChatGPT.
Also ChatGPT was used for looking for typos in the formulas in the PDF file and the converted LaTeX file.
WinEdt was used for spell checking.

\bigskip
\bigskip

\begin{center}
\textbf{References}
\end{center}

\medskip

\begin{list}{}{%
  \setlength{\leftmargin}{3.5em}
  \setlength{\labelwidth}{3em}
  \setlength{\labelsep}{0.5em}
  \setlength{\itemsep}{4pt}
  \setlength{\parsep}{0pt}
}
\item[\mbox{[Bo1]}] Borovoi, M. V.: The algebraic fundamental group and abelian Galois cohomology of reductive algebraic groups. Preprint MPI/89--90, Bonn, 1989.

\item[\mbox{[Bo2]}] Borovoi, M. V.: Abelianization of the second non-abelian Galois cohomology. Duke Math. J. 72 (1993), no. 1, 217--239.

\item[\mbox{[Bo3]}] Borovoi, M. V.: Abelian Galois cohomology of reductive groups. Mem. Amer. Math. Soc. 132 (1998), no. 626, viii+50 pp.

\item[\mbox{[Br1]}] Breen, L.: Bitorseurs et cohomologie non ab\'{e}lienne. In: Cartier, P. et al (eds.) The Grothendieck Festschrift, vol. 1, pp. 401--476, Boston: Birkha\"{u}ser 1990.

\item[\mbox{[Br2]}] Breen, L.: Th\'{e}orie de Schreier sup\'{e}rieure. Ann. Sci. \'Ecole Norm. Sup. (4) 25 (1992), no. 5, 465--514.

\item[\mbox{[Br3]}] Breen, L.: Letter to P. Deligne of February 4, 1988.

\item[\mbox{[BHi]}] Brown, R., Higgins, P. J.: Crossed complexes and non-Abelian extensions, In: Kamps, K. H. et al (eds.) Category Theory (Lecture Notes in Math. vol. 962, pp. 39--50). Berlin: Springer 1982.

\item[\mbox{[BHu]}] R. Brown. R, Huebschmann, P.: Identities among relations, In: Brown, R., Thickstun, T. L. (eds.) Low-Dimensional Topology, (London Math. Soc. Lecture Notes Series vol. 48, pp. 153--202). Cambridge: Cambridge University Press 1982.

\item[\mbox{[Co]}] D. Conduch\'{e}, Modules crois\'{e}s g\'{e}n\'{e}ralis\'{e}s de longueur 2. J. of Pure and Applied Algebra \textbf{34}, 155--178 (1984).

\item[\mbox{[Ded1]}] Dedecker, P.: Sur la cohomologie non ab\'{e}lienne I. Can. J. Math. \textbf{12}, 231--251 (1960).

\item[\mbox{[Ded2]}] Dedecker, P.: Les foncteurs $\mathrm{Ext}_{\Pi}$, $\mathbf{H}_{\Pi}^{2}$ et $H_{\Pi}^{2}$ non ab\'{e}liens. C. R. Acad. Sci. \textbf{258}, 4891--4894 (1964).

\item[\mbox{[Ded3]}] Dedecker, P.: Three dimensional non-abelian cohomology for groups. In: Hilton, P. (ed.) Category Theory, Homology Theory and Their Applications II. (Lecture Notes in Math. vol. 92, pp. 32--64). Berlin: Springer 1969.
\end{list}

\begin{list}{}{%
  \setlength{\leftmargin}{3.5em}
  \setlength{\labelwidth}{3em}
  \setlength{\labelsep}{0.5em}
  \setlength{\itemsep}{4pt}
  \setlength{\parsep}{0pt}
}
\item[\mbox{[Del]}] Deligne, P.: Vari\'{e}t\'{e}s de Shimura: Interpr\'{e}tation modulaire, et techniques de construction de mod\`{e}les canoniques. In: Borel, A., Casselman, W. (eds.) Automorphic Forms, Representations, and $L$-functions, (Proc. Symp. Pure Math. vol. 33, Part 2, pp. 247--289). Providence: AMS 1979.

\item[\mbox{[Gi]}] Giraud, J.: Cohomologie non ab\'{e}lienne. (Grundlehren der mathematischen Wissenschaften vol. 179). Berlin: Springer 1971.

\item[\mbox{[Kn]}] Kneser, M.: Galoiskohomologie halbeinfacher algebraischer Gruppen \"{u}ber $\mathfrak{p}$-adischen K\"{o}rpern I. Math. Z. \textbf{88} 40--47 (1965); II. Math. Z. \textbf{89} 250--272 (1965).

\item[\mbox{[Ko1]}] Kottwitz, R. E.: Stable trace formula: cuspidal tempered terms. Duke Math. J. \textbf{51}, 611--650 (1984).

\item[\mbox{[Ko2]}] Kottwitz, R. E.: Stable trace formula: elliptic singular terms. Math. Ann. \textbf{275}, 365--399 (1986).

\item[\mbox{[Mi]}] Milne, J. S.: The points on a Shimura variety modulo a prime of good reduction. In: Langlands, R. P., Ramakrishnan, D. (eds.) The Zeta Functions of Picard Modular Surfaces, pp. 151--253. Montr\'{e}al: Les Publications CMS, 1992.

\item[\mbox{[Noo]}] Noohi, B.: Group cohomology with coefficients in a crossed module.
J. Inst. Math. Jussieu {\bf 10}, 359--404 (2011).

\item[\mbox{[Se]}] Serre, J.-P.: Cohomologie galoisienne. (Lecture Notes in Math. vol. 5). Berlin: Springer 1964.

\item[\mbox{[Sp]}] Springer, T. A.: Non-abelian $H^{2}$ in Galois cohomology. In: Borel, A., Mostow, G. D. (eds.) Algebraic Groups and Discontinuous Subgroups. (Proc. Symp. Pure Math. vol. 9, pp. 164--182). Providence: AMS 1966.

\item[\mbox{[W1]}] Whitehead, J. H. C.: Note on a previous paper entitled ``On adding relations to homotopy groups''. Ann. of Math. \textbf{47}, 806--810 (1946).

\item[\mbox{[W2]}] Whitehead, J. H. C.: Combinatorial homotopy II. Bull. Amer. Math. Soc. \textbf{55}, 453--496 (1949).
\end{list}

\vfill

\end{document}